\documentclass[12pt, A4]{article}
\usepackage[utf8]{inputenc}

\usepackage{url}
\RequirePackage{filecontents}
\usepackage{blindtext}
\usepackage{amssymb}
\usepackage{amsthm}
\usepackage{amsmath}
\usepackage{mathtools}
\usepackage{xcolor}
\usepackage{sectsty}
\usepackage{hyperref}

\sectionfont{\fontsize{15}{15}\selectfont}

\usepackage{geometry}
\numberwithin{equation}{subsection}
\renewcommand*{\theequation}{%
  \ifnum\value{subsection}=0 %
    \thesection
  \else
    \thesubsection
  \fi
  .\arabic{equation}%
}

\swapnumbers

\theoremstyle{plain}
\newtheorem{thm}[subsection]{Theorem}
\newtheorem{conj}[subsection]{Conjecture}

\theoremstyle{definition}

\newcommand{\R}{\mathbb{R}} 
\newcommand{\N}{\mathbb{N}} 
\newcommand{\A}{\mathcal{A}} 
\newcommand{\B}{\mathcal{B}} 
\newcommand{\AB}{\mathcal{A}\otimes\mathcal{B}} 
\newcommand{\sigalg}{$\sigma$-algebra } 

\newcommand{\brac}[1]{\left( #1 \right)} 

\makeatletter
\def\moverlay{\mathpalette\mov@rlay}
\def\mov@rlay#1#2{\leavevmode\vtop{%
   \baselineskip\z@skip \lineskiplimit-\maxdimen
   \ialign{\hfil$\m@th#1##$\hfil\cr#2\crcr}}}
\newcommand{\charfusion}[3][\mathord]{
    #1{\ifx#1\mathop\vphantom{#2}\fi
        \mathpalette\mov@rlay{#2\cr#3}
      }
    \ifx#1\mathop\expandafter\displaylimits\fi}
\makeatother

\newcommand{\bigcupdot}{\charfusion[\mathop]{\bigcup}{\cdot}}

\title{The translation invariant product measure problem in non-sigma finite case}
\author{Nicha Khenkhok \\ \small Faculty of Science, Eötvös Loránd University \\\small Pázmány Péter stny. 1/A, 1117 Budapest, Hungary;\\ \small e-mail: knicha@student.elte.hu}

\begin{document}

\maketitle

\begin{abstract}
    We give an example of non-translation invariant product measure obtained from two translation invariant measures, one of which is non-sigma finite. This particular example also suggests that there can be infinitely many product measures if we abandon the sigma-finiteness assumption. 
\end{abstract}

\section{Introduction}
Let $\brac{X,\A,\mu}$ be a measure space with \sigalg $\A$ and measure $\mu$, and $\brac{\R,\B,\nu}$ be a real measure space with the Borel \sigalg $\B$ and the Lebesgue measure $\nu$. Denote their product measure space by $\brac{X\times\R, \AB, \mu\times\nu}$, where the product measure is arbitrary. Define a product measure using the definition given by D.H. Fremlin in \cite{fremlin}. The set function $\mu\times\nu:\AB\to[0,\infty]$ is a product measure iff it is a measure and for every measurable rectangle $A\times B$, where $A\in\A$ and $B\in\B$, we have
\[
    \mu\times\nu\left(A\times B\right) = \mu(A)\nu(B).
\]
We shall fix these measure spaces throughout the article.\\

The Lebesgue measure is known to be translation-invariant. One question we may ask is whether a product measure $\mu\times\nu$ inherits this property in the sense that any shift applied to a measurable set $B\in \AB$ along the real axis does not alter the measure. Formally, we conjecture
\begin{conj}
Let the product measure space $\brac{X\times\R, \AB, \mu\times\nu}$ be arbitrary and a set $B\in\AB$ be given. For any $c\in \R$, define the vertical shift of $B$ by $c$ as the set
$$B+c \coloneqq \left\{(x,y+c): (x,y)\in B\right\}\in \AB.$$
Then, $\mu\times\nu\left(B+c\right) = \mu\times\nu\left(B\right)$.
\end{conj}

If the measure space $\brac{X,\A,\mu}$ is $\sigma$-finite, then the conjecture holds trivially as the product measure is unique. This unique product measure is obtained through the Carathéodory's extension theorem. As for the non-$\sigma$-finite case, we will show that the conjecture is not true. 

\section{Completely locally determined product measure}

Let $(X,\A,\mu)=([0,1],\B,\mu)$, where $\B$ is the Borel $\sigma$-algebra and $\mu$ is the counting measure. Then, we may define the measurable space of $(X,\A,\mu)$ and $(\R,\B,\nu)$. Let 
\[
\pi(E) = \inf\left\{\sum_{n=0}^\infty \mu\times\nu\left(A_n\times B_n\right ): \{A_n\}_{n\in\N}\subseteq X, \{B_n\}_{n\in\N}\subseteq \R, E\subseteq \bigcup_{n=0}^\infty A_n\times B_n \right\}
\]
be the product measure space obtained through the Carathéodory's extension theorem.\\

Another candidate as a product measure is the completely locally determined product measure (c.l.d), which the reader may refer to \cite{fremlin} for further details. The c.l.d product measure is given by 
\[
\rho(E)=\left\{\pi(E\cap(A\times B)): A\in \A, B\in \B, \mu(A)<\infty, \nu(B)<\infty\right\}.
\]

On the diagonal $\Delta=\{(x,x):x\in [0,1]\}$, which can be written as
\[
    \Delta = \bigcap_{n=1}^\infty\bigcup_{k=0}^\infty \left[\frac{k}{n},\frac{k+1}{n}\right]\times\left[\frac{k}{n},\frac{k+1}{n}\right] \in \AB,
\]
we have $\pi(\Delta)=\infty$ and $\rho(\Delta)=0$.

 \section{Counterexample measure}
 We will construct a product measure, which utilises the c.l.d. measure. Let $\Delta=\{(x,x): x\in[0,1]\}$ as before. Recall that $\nu:\B\to[0,\infty]$ is the Lebesgue measure on the Borel $\sigma$-algebra. Define $f:[0,1]\to[0,1]\times[0,1]$ to be 
 \[
 f(x)=(x,x),
 \]
 which is a measurable function on $[0,1]$. As every preimage $f^{-1}[E]$ of a measurable set $E\in\AB$ is measurable in $\B$, we can safely define the set function $\xi:\AB\to[0,1]$ as
 \[
 \xi(E) = \nu(f^{-1}[E\cap\Delta]).
 \]
 We claim that $\xi$ is a measure. Trivially, $\xi(\emptyset)=0$. We now check the $\sigma$-additivity property. Let $\{E_n\}_{n\in\N}\subseteq\AB$ be a sequence of disjoint sets. Then, 
 \begin{align*}
     \xi\left(\bigcupdot_{n=0}^\infty E_n\right)
     &=\nu\left(f^{-1} \left[\bigcupdot_{n=0}^\infty E_n\cap \Delta\right]\right)\\
     &= \nu\left(f^{-1} \left[\bigcupdot_{n=0}^\infty (E_n\cap \Delta)\right]\right)\\
     &= \nu\left( \bigcupdot_{n=0}^\infty f^{-1}[E_n\cap \Delta]\right) \\
     &= \sum_{n=0}^\infty \nu\left(f^{-1}[E_n\cap \Delta]\right)\\
     &= \sum_{n=0}^\infty \xi(E_n).
 \end{align*}
 That is, $\xi$ is indeed a measure on $\AB$. We now proceed to the main result.
 \begin{thm}
    There exists a product measurable space $\brac{X\times\R, \AB, \mu\times\nu}$ such that for some $c\in\R$ and some measurable set $B\in\AB$, the vertical shift of $B$ by $c$ results in a change in measure. That is, $\mu\times\nu(B)\neq\mu\times\nu(B+c)$.
 \end{thm}
 \proof 
 We assume the notions previously defined in this section.  Consider the set function $\eta:\AB\to[0,\infty]$ given by
 \[
 \eta(E) = \rho(E)+\xi(E).
 \]
 Since $\eta$ is a sum of measures on $\AB$, we have that $\eta$ is also a measure on $\AB$. We remain to prove that $\eta$ is a product measure. For this, we consider the following cases for a measurable rectangle $A\times B$, where $A\in\A$ and $B\in\B$.
 \begin{itemize}
     \item If $\mu(A)<\infty$ and $\nu(B)\leq \infty$, then $A$ has finitely many points since $\mu$ is a counting measure. So, $A= \{a_1,...,a_k\}$ for some $k\in\{0,1,\dots\}$. It holds that
     \[
     A\times B = \{a_1,...,a_k\}\times B \subseteq \{a_1,...,a_k\}\times \R = A\times \R ,
     \]
     and hence,
     \[
        \Delta\cap \left( A\times B\right)\subseteq \Delta\cap \left( A\times \R\right) = \{(x,x): x=a_1,...,a_k\}.
     \]
     Using monotonicity of measure, 
     \[
     \xi(A\times B) \leq \xi(A\times \R ) = \nu(f^{-1}[\Delta\cap \left( A\times \R\right)])=\nu(\{a_1,...,a_k\})=0.
     \]
     Therefore, $\eta(A\times B) = \rho(A\times B ) + \underbrace{\xi(A\times B)}_{0} = \rho(A\times B) = \mu(A)\nu(B)$.
     \item If $\mu(A)=\infty$ and $\nu(B)> 0$, then $\rho(A\times B) = \mu(A)\nu(B) = \infty$. Therefore, 
     \[
     \eta(A\times B) = \underbrace{\rho(A\times B)}_{\infty} + \underbrace{\xi(A\times B)}_{\geq 0} = \underbrace{\rho(A\times B)}_{\infty}  = \mu(A)\nu(B).
     \]
     \item If $\mu(A)=\infty$ and $\nu(B)= 0$, then $\rho(A\times B) = \mu(A)\nu(B) = 0$. It holds that
     \[
     f^{-1}[\Delta\cap(A\times B)]\subseteq f^{-1}[\Delta\cap(\R \times B)]=B\cap[0,1]
     \]
     By monotonicity of measure, 
     \[
        \xi(A\times B)=\nu(f^{-1}[\Delta\cap(A\times B)])\leq \nu(B\cap[0,1])\leq \nu(B)= 0.
     \]
     Thus, $\eta(A\times B)=\rho(A\times B) +\xi(A\times B) = 0 = \mu(A)\nu(B)$.
 \end{itemize}
 Therefore, $\eta$ is indeed a product measure. Furthermore, $\eta(\Delta)= \rho(\Delta)+\xi(\Delta) = 0 + 1 = 1$. However, $\eta(\Delta+1)= \rho(\Delta+1)+\xi(\Delta+1)= 0+0 =0$. $\blacksquare$

\end{document}